\documentclass[twoside, 11pt]{article}
\usepackage{mathrsfs}

\usepackage{mathrsfs,amsfonts,amsmath}
\RequirePackage{ifthen,calc}
\usepackage{theorem}
\usepackage[dvips]{color}

 \catcode`@=11
 \def\@evenhead{\hbox to\textwidth{\footnotesize\rm\thepage \hfill
  {\it }}} 

 \def\@oddhead{\hbox to \textwidth{\footnotesize{\it
  ~  } \hfill\thepage}}

 \renewcommand{\section}{\makeatletter
 \renewcommand{\@seccntformat}[1]{{\csname the##1\endcsname.}\hspace{0.45em}}
 \makeatother \@startsection
{section}
{1}
{0pt}
{\baselineskip}
{0.5\baselineskip}
{\normalsize\bfseries\mathversion{bold}}}

\catcode`@=12
\newcommand\ack{\section*{Acknowledgement}}

\newtheorem{thm}{\noindent Theorem}[section]
\newtheorem{lem}{\noindent Lemma}[section]
\newtheorem{cor}{\noindent Corollary}[section]

\newtheorem{remark}{\noindent Remark}[section]

\newtheorem{proposition}{Proposition}[section]
\newtheorem{lemma}[proposition]{Lemma}
\newtheorem{theorem}[proposition]{Theorem}

\newtheorem{corollary}[proposition]{Corollary}
\newtheorem{example}[proposition]{Example}
\theoremheaderfont{\normalfont\bfseries} \theorembodyfont{\slshape}
{\theorembodyfont{\rmfamily}}
{\theorembodyfont{\rmfamily}}
{\theorembodyfont{\rmfamily}}
{\theorembodyfont{\rmfamily}}
{\theorembodyfont{\rmfamily}}
{\theorembodyfont{\rmfamily}}
{\theorembodyfont{\rmfamily}}
{\theorembodyfont{\rmfamily}}
{\theorembodyfont{\rmfamily}}

 \def\beqlb{\begin{eqnarray}}\def\eeqlb{\end{eqnarray}}
 \def\beqnn{\begin{eqnarray*}}\def\eeqnn{\end{eqnarray*}}
 \newcommand{\bgeqn}{\begin{equation}}
\newcommand{\edeqn}{\end{equation}}
\def\ra{\rightarrow}

 \numberwithin{equation}{section}
\def\qed{\hfill$\square$\smallskip}
\def\no{\nonumber}
\def\L{{\mathcal L}}
\def\R{{\mathbb R}}
\def\N{{\mathbb N}}
\def\Q{{\mathbb Q}}
\def\Z{{\mathbb Z}}
\def\M{{\EuScript M}}
\def\EN{{\EuScript{E}}}
\def\CAP{{\EuScript{C}}}
\def\ra{\rightarrow}
\def\iy{\infty}
\def\sleq{\mathop{\preccurlyeq}}
\def\sgeq{\mathop{\succcurlyeq}}
\def\bfE{{\mathbb{E}}}
\def\mbe{{\mathbb{E}}}
\def\mbfE{{\mathbf{E}}}
\def\mbr{{\mathbb{R}}}
\def\bfP{{\mathbb{P}}}
\def\mbP{{\mathbb{P}}}
\def\mbfP{{\mathbf{P}}}
\def\bfR{{\mathbb{R}}}
\def\bfQ{{\mathbb{Q}}}
\def\bfN{{\mathbb{N}}}
\def\1{{\mathbf{1}}}
\def\e{\mathrm{e}}
\def\d{\mathrm{d}}
\def\red{\color{red}}
\def\blue{\color{blue}}
\def\ll{{\langle}}
\def\r{\rangle}
\def\eps{\varepsilon}
\def\Re {{\rm Re}\,}
\def\mcr{\mathscr}
\def\mtf{{\mathcal{F}}}
\def\pos{\bar{{\mbb R}}_{+}}
\def\pint{\int_{0}^{+\infty}}
\def\nint{\int_0^{+\infty}}
\def\myto{\stackrel{\mathcal{D}}{\to}}

\begin{document}
\title{\bf Weak Quenched Invariance Principle for Random Walk with Random Environment in Time
}
\author{You Lv\thanks{Email: lvyou@dhu.edu.cn },~~Wenming Hong\thanks{Corresponding author. Email: wmhong@bnu.edu.cn}
\\
\\ \small College of Science, Donghua University,
\\ \small Shanghai 201620, P. R. China.
\\
\\ \small School of Mathematical Sciences $\&$ Laboratory of Mathematics and Complex Systems,
\\ \small Beijing Normal University, Beijing 100875, China
}
\date{}
\maketitle

\noindent\textbf{Abstract} Consider the invariance principle for a random walk with random environment (denoted by $\mu$) in time on $\bfR$ in a weak quenched sense. We show that a sequence of the random probability measures on $\bfR$ generated by a bounded Lipschitz functional $f$ and $\mu$ will converge in distribution to another random probability measures, which is related to $f$ and two independent Brownian motions. The upper bound of the convergence rate has been obtained. We also explain that in general, this convergence can not be strengthened to the almost surely sense.

\smallskip


\smallskip

\section{Introduction}
Donsker's classical invariance principle (see \cite{D1951}) tells that the random broken line connecting vertexes $\left(i/n, Y_i/\sqrt{n}\right), i\leq n$ will converge in distribution to a standard Brownian motion $\{B_s\}_{s\in[0,1]}$
~if $\{Y_n\}_{n\in\bfN}$ is an i.i.d. random walk with $\bfE(Y_1)=0$ and $\bfE(Y^2_1)=1.$ In the present paper, we extend this result to the model called random walk with random environment in time. Now we first give the definition of this model.

Denote $\mu:=\{\mu_n\}_{n\in\bfN^+}$ an i.i.d. sequence with values in the space of probability measures on $\bfR.$ Conditioned to a realization 
of $\mu,$ we sample $\{X_n\}_{n\in\bfN^+}$ a sequence of independent random variables such that for every $n,$ the law of $X_n$ is the realization of $\mu_n.$
Set \beqlb\label{sec-sn} S_0=0\in\bfR,~~~S_n:=S_0+\sum_{i=1}^{n}X_i.\eeqlb 
We call $\{S_n\}_{n\in\bfN}$ the {\it random walk with random environment $(\mu)$ in time}.  
There are two laws to be considered for this model. 
We write $\mbfP_{\mu}$ for the law of  $\{S_n\}_{n\in\bfN}$ conditionally on the sequence $\{\mu_n\}_{n\in\bfN^+},$ which is usually called a quenched law. Let $\mathbf{P}$ be the joint law of $\{S_n\}_{n\in\bfN}$ and the environment $\mu,$ which is usually called an annealed law.
The corresponding expectations are denoted by $\mbfE_{\mu}$ and $\mathbf{E}$ respectively. In the quenched sense, $\{S_n\}_{n\in\bfN}$ can be viewed as a sequence of sums of independent, not necessarily identically distributed random variables. Obviously, if the random environment is degenerate, then $\{S_n\}_{n\in\bfN}$ will degenerate to an i.i.d. random walk. To distinguish from the notations in the case of random environment, throughout the present paper we write $\bfP, \bfE$ the probability and the corresponding expectation for the random walk with fixed environment (including time-homogeneous and time-inhomogeneous case).

Note that the process considered in the present paper is not the
random walk in random environment which has been well-studied in Zeitouni \cite{Z2004}, Goldsheid \cite{GI} and many other papers, where the i.i.d. random environment (some papers consider the stationary ergodic random environment) is either purely spatial or space-time. Contrary to this setting, in our model the random environment varies in time randomly. In the rest of the present paper,  we often abbreviate the model we consider as RWre rather than RWRE, since the latter abbreviation has been used to present the model in Zeitouni \cite{Z2004}.

A connection between this two models can be stated as follows. For a one-dimension RWRE $\{Z_n\}$ with only nearest neighbour jumping (i.e., for almost all environments and any $n$, the quenched probability of $\{|Z_{n+1}-Z_{n}|\neq 1\}$ is $0$), the hitting time
\beqlb\label{htime}T_n:=\inf\{i: Z_i=n\}\eeqlb
is a RWre when the random environment of $Z_n$ is i.i.d.. But we can see $\{T_n\}$ is a special one since it has strict increasing trajectories for all environments and $\tau_n:=T_{n+1}-T_n$ can only take odd numbers since the nearest neighbour jumping, which is an essential setting for RWRE in many papers. While the RWre we consider is more generalized than $\{T_n\},$ see Section 2 for details. We compare our result with some achieved conclusions on $T_n$ and give a discussion on the moment condition after Corollary \ref{t02}.

The motivation to study RWre comes from the research on a series of limit behaviors of the branching random walk with a time-inhomogeneous random environment (BRWre) on $\bfR$, which is a generalized branching random walk. Compared with the time-homogeneous branching random walk, BRWre stress that the point processes of the reproduction law (including the branching and displacement) of the particle(s) in each generation are random (according to a given distribution on the collection of all point processes on $\bfR$ independently) but not a fixed one. The time-inhomogeneous version of Many-to-one formula given in Mallein \cite{M2015a} constructs a bijection between RWre and BRWre, which makes some problems on the asymptotic behaviors of BRWre is equivalent to the problems on some asymptotic behaviors of the quenched distributions of RWre. Hence RWre can be regarded as an essential tool to study BRWre.
For example, to obtain the second order of the maximal displacement of a BRWre (see Mallein and Mi{\l}o\'{s} \cite{MM2016}), Mallein and Mi{\l}o\'{s} \cite{MM2015} introduced the RWre $\{S_n\}$ and show that
\beqlb\label{1.2RBT}\lim\limits_{n\ra \infty}\frac{\log\mbfP_{\mu}(\forall_{i\leq n} S_i\geq 0|S_0=x)}{\log n}=-\gamma_1,~~ \forall x>0, ~~\rm \mathbf{P}\text{-}a.s.,\eeqlb where $\gamma_1$ and the coming $\gamma_2$ are constants depending on the distribution of $S_1.$ To investigate the barrier problem of a BRWre (see Lv \& Hong \cite{LY3}), a previous work of the authors \cite{LY2} gave the quenched small deviation for the trajectories of a RWre, where we showed that
\beqlb\label{1.2SD}\lim\limits_{n\ra \infty}\frac{\log\mbfP_{\mu}(\forall_{i\leq n} |S_i|\leq n^{\alpha}|S_0=x)}{n^{1-2\alpha}}=-\gamma_2,~~ \forall x\in\bfR, \alpha\in\big(0,\frac{1}{2}\big), ~~\rm \mathbf{P}\text{-}a.s.,\eeqlb
which extends the main result (small deviation for random walk) in Mogul'ski\u{\i} \cite{Mog1974}.

In fact, using the time-homogeneous version of Many-to-one formula and some conclusions on the distributions of time-homogeneous random walk to study BRW is a standard and effective way which has been proven in many papers, see \cite{AJ2011}, \cite{GHS2011} and \cite{HS2009}. Hence we believe that the more investigations on RWre will contribute to the future exploration on BRWre. The present paper focuses on the quenched invariance principle of the RWre, which is inspired by the following observation of \eqref{1.2RBT} and \eqref{1.2SD}. We think that the reason why the quenched probabilities in \eqref{1.2RBT} and in \eqref{1.2SD} converge to $0$ and have a $\mathbf{P}\text{-a.s.}$ decay rate is the time-space scaling therein. Moreover, we see the probability in \eqref{1.2SD} will decay more slowly as $\alpha\ra \frac{1}{2}.$ Then a natural question is, what will happen if we set a proper time-space scaling for RWre? In the light of \eqref{1.2SD}, we consider the scaling as the case $\alpha=\frac{1}{2}.$ We begin by introducing a heuristic example as follows.



Consider a special type of RWre which we tend to call it bio-normal RWre. Assume that the quench mean of $S_1$ is of a normal distribution $\mathcal N(0, \sigma_1^2)$; moreover, for any realization of $\mu,$ $S_1-\mbfE_{\mu}(S_1)$ has a common normal distribution $\mathcal N(0, \sigma_2^2)$ under $\mbfP_{\mu}$ ($\sigma_1$ and $\sigma_2$ are both constants). Note that
 $$\mbfP_{\mu}\left(\max_{i\leq n}|n^{-1/2}S_i|\leq 1\right)=\mbfP_{\mu}\left(\forall_{i\leq n}S_i-\mbfE_{\mu}S_1\in [-n^{1/2}-\mbfE_{\mu}S_1, n^{1/2}-\mbfE_{\mu}S_1]\right),$$
 then by the scaling property of Brownian motion, we see $\mbfP_{\mu}\left(\max_{i\leq n}|n^{-1/2}S_i|\leq 1\right)$ has the same distribution as
 $\bfP\left(\max_{i\leq n}|\sigma_2B_{i/n}+\sigma_1W_{i/n}|\leq 1|W\right),$ where $B$ and $W$ are two independent standard Brownian motions.
From this view we wonder whether for a more general RWre, $\mbfP_{\mu}\left(\max_{i\leq n}|n^{-1/2}S_i|\leq 1\right)$ will converge to $\bfP\left(\sup_{s\leq 1}|r_2B_s+r_1W_s|\leq 1|W\right)$ in distribution (under some mild assumption), where $r_1, r_2$ are two constants. We give the answer in remark \ref{rm2.1r}.
 The main result of the present paper is to give a more general form of the invariance principle by the bounded Lipschitz functional on $C[0,1].$ (It is well known that the weak convergence in a metric space can be characterized equivalently in several ways including the form of Lipschitz functional, see \cite[Theorem 3.10.1]{Dbook}.) Denote $S_{(n)}$ the random broken line connecting vertexes $\{\left(i/n, S_i/\sqrt{n}\right)\}_{i\leq n}.$ Note that for any bounded Lipschitz functional $f$ defined on $C[0,1],$ the distribution of $\mbfE_{\mu}(f(S_{(n)}))$ can be viewed as a random probability measure on $\bfR.$
 We show that there exists two constants $r_1,r_2$ such that $\mbfE_{\mu}(f(S_{(n)}))$ will converges to $\mbfE(f(r_1B+r_2W)|W)$ in distribution, where $r_1B+r_2W$ and $W$ present the random trajectories $\{r_1B_s+r_2W_s, s\in[0,1]\}$ and $\{W_s, s\in[0,1]\}$ respectively. Moreover, an upper bound of the convergence rate is obtained.

 Refining the idea from \cite{MM2015} and \cite{LY2}, the main tool we use in the present paper is the strong approximation of sums of independent variables with finite moments (see Theorem \uppercase\expandafter{\romannumeral1}, Section 3.1). According to the strong approximation theorem, we can divide the environment space into the ``regular" part (in which we can find a Brownian motion to approximate the sums of independent variables) and the ``irregular" parts (which happens with very small probability). The difficulties in this work are how to estimate the random functional of the random environment and how to deal with the the interdependence between $\mbfE_{\mu}S_n$ and $S_n-\mbfE_{\mu}S_n$. We also point out that our approach is different from the ones (e.g. heat kernel estimates, pointwise sublinearity of corrector, Lindeberg-Feller martingale functional CLT) used in the references on the invariance principle for RWRE (see \cite{BP2007, NK2013, RS2009}) or some other random motions in random media (see \cite{BB2007,SS2004,SZ2006}).

 It should be stressed that the diffusion constant/matrix of the Brownian motions as the limits of the convergence in the references mentioned above are independent of the realization of the random media. In this paper, we will see that the limit of the quenched invariance principle will be a  Brownian motion with a random center and the center depends on the realization of the random environment.

 In the next section, we give the main theorem, corollaries and remarks. The proof is given in Section 3.

\section{Main result}
The RWre $\{S_n\}_{n\in\bfN}$ (with the random environment $\mu$) is defined as \eqref{sec-sn}. Denote \begin{eqnarray}\label{basic}M_{n}:=\mbfE_\mu(S_n), ~U_{n}:=S_n-\mbfE_\mu(S_n),~\eta_n:=\mbfE_\mu(U^2_n)=\mbfE_\mu(S^{2}_n)-M^2_{n}.\end{eqnarray}
Throughout this paper we assume that $S_0=0, \mbfE(M_{1})=0, \mbfE(\eta_1)=1$ \footnote{Note that $\mbfE(\eta_1)=0$ means that the RWre degenerates to i.i.d. random walk, since in the sense of ``up to a constant'', the assumption $\mbfE(\eta_1)=1$ is without loss of generality.}. Denote $\sigma^2:=\mbfE(M^2_{1}).$
Obviously, it is true that $\mbfE(U_n^2)=\mbfE(\eta_n).$

We set $f$ a bounded Lipschitz function on $C[0,1]$ with respect to the uniform norm, i.e., there exist constants $K, L$ such that $$|f(x)|\leq L~~\text{and}~~|f(x)-f(y)|\leq K\sup_{t\in[0,1]}|x(t)-y(t)|,~~ \forall x,y\in C[0,1].$$
Throughout the present paper, we denote $S_{(n)}\in C[0,1]$ be the random broken line generated by the random walk $\{S_n\}$ in the way
$$S_{(n)}(t):=\frac{S_{\lfloor tn\rfloor}+(tn-\lfloor tn\rfloor)(S_{\lfloor tn\rfloor+1}-S_{\lfloor tn\rfloor})}{\sqrt{n}},~~t\in[0,1],$$
where $\lfloor \cdot\rfloor$ represents the maximal integer value not exceeding $\cdot.$
For a Brownian motion $Y,$ slightly abusing notation we usually write $\{Y_t, t\in[0,1]\}$ as $Y$ for simplicity.
For $A, A'\in C[0,1], c\in\bfR,$ $A+A'$ presents $(A+A')(t):=A(t)+A'(t)$ and $cA$ presents $(cA)(t):=cA(t),\forall t\in[0,1].$
Denote $F_{n}$ and $F$ the distribution functions of $\mbfE_{\mu}(f(S_{(n)}))$ and $\mbfE(f(B+\sigma W)|W),$ where $B$ and $W$ are two independent standard Brownian motions on $\bfR$ with $B_0=W_0=1$. Now we state the main result of this paper.
\begin{theorem}\label{t1}
Assume that there exist $\alpha,\beta_1>2, \beta_2\geq 2, \beta_3\geq 2$ such that \begin{equation}\label{ast1}\mbfE(|M_1|^{\beta_1})+\mbfE(\eta_1^{\beta_2})+\mbfE(\zeta_1^{\beta_3})<+\infty,\end{equation}
where $\zeta_n:=\sum_{k=1}^n\mbfE_\mu(|U_k-U_{k-1}|^{\alpha}).$ Sequences $\{a_n\},\{b_n\}$ satisfy that
\begin{equation}\label{a2a}  \tilde{a}:=\varliminf_{n\ra\infty}\frac{a_n}{n^{-1/4}(\log n)^{3/4}}>0,~\lim\limits_{n\ra\infty}a_n=0;\end{equation}
\begin{eqnarray}\label{ab}\forall n\in\bfN,~ b_n\geq 2\mbfE\zeta_1n , ~~\varlimsup_{n\ra\infty}b_n(a_n\sqrt{n})^{-\alpha}=0.~\end{eqnarray}
Then there exists two positive constants $K_1$ (depending on both $f$ and the distribution of $S_1$, specifically, the values of $K, L, \alpha, \beta_2$ and $\tilde{a}$) and $K_2$ (depending only on the distribution of $S_1$, specifically, the values of $\mbfE\left(|M_1|^{\beta_1}\right),~\mbfE\left(|\eta_1-1|^{\beta_2}\right),$ $\mbfE\left(|\zeta_1-\mbfE\zeta_1|^{\beta_3}\right),$ $ \alpha,$ $\beta_1,$ $\beta_2$ and$~\beta_3$) such that
\begin{eqnarray}\label{thu}
\forall n\in\bfN,~x\in\bfR,~~ F_{n}(x)\leq F(x+K_1y_n)+K_2r_n,
\end{eqnarray}
\begin{eqnarray}\label{thl}
\forall n\in\bfN,~x\in\bfR,~~ F_{n}(x)\geq F(x-K_1y_n)-K_2r_n,
\end{eqnarray}
where $y_n:=a_n+b_n(a_n\sqrt{n})^{-\alpha}, r_n:=n^{1-\frac{\beta_1}{2}}a_n^{\beta_1}+n^{1-\beta_2}(\log n)^{\beta_2}a_n^{-2\beta_2}+nb^{-\beta_3}_n.$ 
\end{theorem}


\begin{corollary}\label{t02}
If $\mbfE\left((\mbfE_{\mu}|S^\alpha_1|)^2\right)<+\infty$ for some $\alpha>2,$ then for any bounded Lipschitz functional $f$, we have
$$\mbfE_{\mu}(f(S_{(n)}))\stackrel{d}{\rightarrow}\mbfE(f(B+\sigma W)|W),$$
where $\stackrel{d}{\rightarrow}$ presents the convergence in distribution.
\end{corollary}

\noindent\emph{Proof.} Note that \eqref{ast1} holds when we take $\beta_1=2\alpha, \beta_2=\alpha, \beta_3=\alpha.$ Hence this corollary is a direct result from the observation $\lim_{n\ra +\infty}y_n=\lim_{n\ra +\infty}r_n=0.$\qed
\begin{remark}\label{rm2.1r}The approach used in the proof of Theorem 2.1 can also adapt to show that for any constants $c_1<0<c_2,$
\begin{eqnarray}\label{rm2.1}\mbfP_{\mu}\Big(\forall_{i\leq n}n^{-\frac{1}{2}}S_i\in[c_1,c_2]\Big)\stackrel{d}{\rightarrow}\mbfP(\forall _{s\leq 1}B_s+\sigma W_s\in[c_1,c_2]|W),\end{eqnarray}
though the indicator function is not a Lipschitz function. Comparing \eqref{rm2.1} with small deviation \eqref{1.2SD}, one see that $\alpha=\frac{1}{2}$ is the proper scaling.
\end{remark}

We point out that the assumption \eqref{ast1} we require is stronger than $\mbfE S^2_1<+\infty$ since we consider the invariance principle of $\{S_n\}$ but not $\{S_n-\mbfE_{\mu}S_n\}$ and we apply the strong approximation. We need a stronger moment condition to control the random center $\mbfE_{\mu}S_n$ of $S_n$ (note that the center of $S_n-\mbfE_{\mu}S_n$ is $0$) and to overcome the difficulty brought from the interdependence of $S_n-\mbfE_{\mu}S_n$ and $\mbfE_{\mu}S_n.$ We may expect a weaker assumption if we only consider $\{S_n-\mbfE_{\mu}S_n\}.$ For example,  Alili \cite{A1999} gave the central limit theorem for the sequence of hitting times $\{T_n\}$ defined in \eqref{htime}. As we have mentioned in Section 1, $\{T_n\}$ is a special kind of RWre (strictly increasing, the increment takes values on the collection of positive odd numbers). Under the conditions \footnote{The condition $\mbfE(T_1)<+\infty$ is not given directly in \cite{A1999} but from the branching structure of the RWRE we can see it is equivalent to the condition ``$E(S)<+\infty$" in \cite[Section 5]{A1999}.} \beqlb\label{conht}\iota^2:=\mbfE[(T_1-\mbfE_{\mu}T_1)^2]<+\infty ~~~\text{and}~~~ \mbfE(T_1)<+\infty,\eeqlb Alili showed that
\begin{eqnarray}\label{rm2.1+}\mbfP_{\mu}\Big((\iota^2 n)^{-\frac{1}{2}}(T_n-\mbfE_{\mu}T_n)\in[c_1,c_2]\Big){\rightarrow}\mbfP(\forall _{s\leq 1}B_s\in[c_1,c_2]|W),~~~\mbfP-{\rm a.s.}.\end{eqnarray}
As far as we know, the central limit theorem \eqref{rm2.1+} is the most relevant one to our results among the achieved conclusions on $\{T_n\}.$
It is indeed that Afanasyev \cite{Af2012} showed an invariance principle for $T_n,$ but it is under another case and hence the limit distribution is different from what we obtain.
In addition, we stress that the method to obtain the central limit theorem, such as checking the Lindeberg condition or the Berry-Esseen estimates can not adapt for the invariance principle given in the present paper. 
Moreover, we give the upper bound of the convergence rate in Theorem \ref{t1} and the following corollary.

\begin{corollary}\label{t2}
Assume that $\mbfE\left((\mbfE_{\mu}|S^\alpha_1|)^2\right)<+\infty$ for some $\alpha>2.$ Let $\pi_n$ be the Prokhorov distance between $\mbfE_{\mu}(f(S_{(n)}))$ and $\mbfE(f(B+\sigma W)|W),$ then we have
$$(1)~ \pi_n\leq O(n^{\frac{1-\alpha/2}{1+\alpha}}) ~~\text{as}~~ \alpha\in(2,5);$$
$$~~(2)~ \pi_n\leq O(n^{-\frac{1}{4}}(\log n)^{\frac{3}{4}}) ~~\text{as}~~ \alpha\geq 5.$$
\end{corollary}
\noindent\emph{Proof.} From Jensen's inequality we see $\mbfE(|M_1|^{2\alpha})<+\infty$ . Note that
$$\mbfE_{\mu}|U^\alpha_1|\leq 2^{\alpha-1}(\mbfE_{\mu}|S^{\alpha}_1|+|M^{\alpha}_1|)<+\infty.$$
Hence $\mbfE((\mbfE_{\mu}|U^\alpha_1|)^{2})<+\infty,$ which means that $\mbfE((\mbfE_{\mu}|U^2_1|)^{\alpha})<+\infty$ by Jensen's inequality. That is to say,
Theorem 2.1 holds when we take $\beta_1=2\alpha, \beta_2=\alpha, \beta_3=2.$ First we consider the case $\alpha\in(2,5).$ Let $a_n:=n^{\frac{1-\alpha/2}{1+\alpha}},$ then $y_n=O(a_n)$ when we take $b_n=O(n).$ Now we only need to show $r_n=o(a_n).$    
Recall the expression of $r_n$ in Theorem 2.1, for this case we see that the order of $r_n$ will be $\max\{n^{1-\beta_2}(\log n)^{\beta_2}a_n^{-2\beta_2},nb^{-\beta_3}_n\}$ since $n^{1-\alpha}a_n^{2\alpha}=o(n^{1-\alpha}(\log n)^{\alpha}a_n^{-2\alpha}).$ Note that in this case $nb^{-2}_n=O(n^{-1})$ hence $\varlimsup_{n\ra\infty}\frac{nb^{-2}_n}{a_n}<+\infty$ and $n^{1-\alpha}(\log n)^{\alpha}a_n^{-2\alpha}=o(n^{1-\frac{\alpha}{2}})$ hence
$n^{1-\alpha}(\log n)^{\alpha}a_n^{-2\alpha}=o(a_n),$ which shows that $r_n=o(a_n)$ and hence $\pi_n\leq O(n^{\frac{1-\alpha/2}{1+\alpha}})$ by the definition of Prokhorov distance. For the case $\alpha\geq 5,$ we take $a_n=n^{-\frac{1}{4}}(\log n)^{\frac{3}{4}}$ and hence $y_n=o(a_n).$ By the above analysis we see $r_n=o(a_n),$ which completes the proof. \qed

\begin{remark}
For an i.i.d. random walk $\{T_n\}$ with $\bfE(T^{\alpha}_1)<+\infty$, $\eqref{sak06}$ (in Section 3) tells that
$$\big|\bfE(f(T_{(n)}))-\bfE(f(B))\big|\leq O\big(n^{\frac{1-\alpha/2}{1+\alpha}}\big),~~ \forall \alpha>2.$$
Comparing this conclusion with Corollary \ref{t2}, we give a comment as follows.
From the proof of Theorem 2.1, we can see the first term in Assumption \eqref{a2a}, which leads that there are two cases in Corollary \ref{t2}, is essentially due to the randomness of $\{\eta_i\}.$ If we assume that $\mbfP(\eta_1=1)=1,$ then we do not need to set the event $J_n$ in Lemma \ref{L1} and hence we do not need to estimate \eqref{thm00+}, which means that we can remove the first restriction on $a_n$ in \eqref{a2a}.
Note that $r_n=o(^{\frac{1-\alpha/2}{1+\alpha}})$ when $\mbfP(\eta_1=1)=1$ (since we can take $\beta_2$ large enough), hence  Corollary \ref{t2} (1) holds for all $\alpha>2.$ That is to say, $\pi_n$ can be very close to $O(n^{-1/2})$ as long as $\alpha$ large enough, which appears reasonable.
This comparison can also be interpreted from a macroscopic view.
We think that RWre has more violent fluctuation than an i.i.d. random walk, since RWre contains two levels of randomness. The common randomness comes from the quenched mean $M_1,$ which plays the role as $T_1$ of the i.i.d. random walk;  while in this comparison, the power of the second randomness is presented by $\eta_1$ when $\alpha\geq 5.$
\end{remark}
\begin{corollary}
If $\mbfE\left((\mbfE_{\mu}|S^\alpha_1|)^2\right)<+\infty$ for some $\alpha>2$ and $\sigma=0,$ then we have $$\mbfE_{\mu}(f(S_{(n)}))\rightarrow\mbfE(f(B)), ~{\rm \mathbf{P}\text{-}a.s.}.$$
\end{corollary}

We remind that $\sigma=0$ does not imply that the random environment is degenerate. It has no impact on the randomness of $\{U_n\}.$

\noindent\emph{Proof.} A key observation is that $\mbfE(f(B+\sigma W)|W)$ is degenerate when $\sigma=0,$ since $B$ is independent of $W.$
If we take $a_n=(\log n)^{-1}, b_n=n\log n,$ then Theorem 2.1 holds for $\beta_1=2\alpha, \beta_2=\alpha, \beta_3=2.$
In this case, one can see $y_n\ra 0$ and $\sum_{n=1}^{+\infty}r_n<+\infty.$ Then we completes the proof by applying Theorem 2.1 and the Borel-Cantelli Lemma.\qed

We use the assumption $\mbfE\left((\mbfE_{\mu}|S^\alpha_1|)^2\right)<+\infty$ rather than $\mbfE\left(|S^{2\alpha}_1|\right)<+\infty$ because the latter is much more restrictive than the former. The following example shows that we can easily construct a RWre satisfying that $\mbfE\left((\mbfE_{\mu}|S^\alpha_1|)^\frac{1}{\varepsilon}\right)<+\infty$ and $\mbfE\left(|S_1^{\alpha+\varepsilon}|\right)=+\infty$ for any $\varepsilon>0.$
\begin{example}
For any realization of $\mu,$ we require $M_1$ is bounded and the distribution function $F_U$ of $U_1$ are the same one, which satisfies that
$F'_U(x)=O(|x|^{-\alpha-1}(\log |x|)^2)$ as $|x|\ra +\infty$.
\end{example}

Most likely, we believe that the convergence in Corollary \ref{t02} can not be strengthened to ${\rm \mathbf{P}\text{-}a.s.}$ if $\mbfE(f(B+\sigma W)|W)$ is non-degenerate. From the proof of Theorem 2.1 and the example ``bio-normal RWre" in Page 2, we see the essential reason can be attributed to the fact that for a Brwonian motion $W,$ $n^{-\frac{1}{2}}W_{sn}$ will not converge to $W_s$ in probability.

\bigskip

\section{Proof}
\subsection{Preliminary}
Let~$\xi_1,\xi_2,\ldots,\xi_n$ be an infinite sequence of independent random variables satisfying $\bfE(\xi_j)=0, \bfE(\xi^2_j)<+\infty$. Denote $D_k:=\sum_{i=1}^k\bfE(\xi^{2}_i).$ Introduce a random broken line $T(s), s\in\bfR^+$ such that $T(0)=0,$ $T(D_k)=\sum_{i=1}^{k}\xi_i, k\in\bfN^+$ and linear, continuous on each interval $[D_{k-1},D_{k}].$ 

\noindent\emph{{\bf Theorem \uppercase\expandafter{\romannumeral1}~(Sakhanenko, \cite[Theorem 1]{Sak2006})}
For any $\alpha\geq 2,$ there exists a standard Brownian motion $B$ (depending on $\alpha$) such that
\begin{eqnarray}\label{sak06}\forall x>0, \bfP\left(\sup_{s\leq D_n}\big|T(s)-B_s\big|\geq x\right)\leq 2(C_{abs}\alpha)^{\alpha}x ^{-\alpha}\sum_{k=1}^n\bfE(\xi^{\alpha}_k),\end{eqnarray}}
where $C_{abs}$ is an absolutely constant.

In the rest of this paper, we denote $C^{(z)}:=2(C_{abs}z)^{z}$ for constant $z$.

\noindent\emph{{\bf Theorem \uppercase\expandafter{\romannumeral2}~(Cs\"{o}rg\H{o}~and~R\'{e}v\'{e}sz, \cite[Lemma 1]{CR1979})}} For a standard Brownian motion $B$ and any $D_1>2$, there exists a constant $D_2>0$ (depending on $D_1$) such that
$$\forall a>0,~t>0,~~~~ \bfP\left(\sup_{0\leq s\leq t}|B_s|\geq a\right)\leq D_2e^{-\frac{a^2}{D_1t}}.$$

We remind that the Berry-Esseen type estimates can not work for our proof.

\subsection{Proof of main theorem}
First we give a lemma which is essential for the proof. Denote $M(s)$ a continuous random broken line which is linear on $[k,k+1]$ and $M(0)=0, M(k)=M_k, k\in\bfN.$
\begin{lemma}\label{L1}
Define the events
$$J_n:=\left\{\forall_{i\leq n}|\eta_{i}-i|\leq \frac{na^2_n}{C_1(\log n)}\right\},~ G_n:=\{\zeta_n\leq b_n\},~$$
where $C_1$ is a constant. Then on the event $J_n\cap G_n,$  these exists a constant $C_2>0$ (depending on $C_1, L, K, \alpha$) such that
\begin{eqnarray}\label{lemu}
\forall n\in\bfN^+,~\Big|\mbfE_{\mu}\left(f(S_{(n)})\right)-\mbfE_{\mu}\left(f(B+M_{(n)})\right)\Big|
\leq C_{2}y_n.\end{eqnarray}
where $y_n:=a_n+b_n(a_n\sqrt{n})^{-\alpha}.$
\end{lemma}
\noindent{\it Proof of Lemma \ref{L1}.}
Recall the notation in \eqref{basic}. Denote $$g(k+\lambda):=(1-\lambda)\eta_k+\lambda\eta_{k+1},k\in\bfN, \lambda\in[0,1).$$ Define a continuous random broken line $U(r), r\in\bfR^+$ by
 $$U(g(s)):=(\lfloor s\rfloor+1-s)U_{\lfloor s\rfloor}+(s-\lfloor s\rfloor)U_{\lfloor s\rfloor+1},s\geq 0.$$
 One can see $U(\cdot)$ is well-defined since $\eta_k=\eta_{k+1}$ implies $\mbfP_{\mu}(U_{k+1}=U_k)=1$.
By Theorem \uppercase\expandafter{\romannumeral1}, we can find a standard Brownian motion $B$ under law $\mbfP_{\mu}$ such that
\begin{eqnarray}\label{del}
\mbfP_{\mu}\left(\max_{s\leq \eta_n}|U(s)-B_s|\geq a_n\sqrt{n}\right)\leq C^{(\alpha)}(a_n\sqrt{n})^{-\alpha}\zeta_n,
\end{eqnarray}
Define $S(s)=(1+\lfloor s\rfloor-s)S_{\lfloor s\rfloor}+(s-\lfloor s\rfloor)S_{\lfloor s\rfloor+1}$ and $M(s)=(1+\lfloor s\rfloor-s)M_{\lfloor s\rfloor}+(s-\lfloor s\rfloor)M_{\lfloor s\rfloor+1},$ ~~$s\in\bfR^+.$
Note that $S(s)=U(g(s))+M(s),$ hence
$$\max_{s\leq n}|S(s)-(B_s+M(s))|\leq \max_{s\leq n}|U(g(s))-B_{g(s)}|+\max_{s\leq n}|B_{g(s)}-B_{\lfloor s \rfloor}|+\max_{s\leq n}|B_{\lfloor s \rfloor}-B_s|.$$
Note that $g(n)=\eta_n,$ hence on the event $J_n$ we have
$|g(s)-s|\leq \frac{na^2_n}{C_1\log n}, \forall s\leq n,$ which means that $|g(s)-\lfloor s \rfloor|\leq \frac{na^2_n}{C_1\log n}+1.$
Choose $C'_0> \sqrt{3\alpha/C_1}.$ Let us further define the event $$\Theta_{n}:=\left\{\forall i\leq n,~\forall |t|\leq \frac{na^2_n}{C_1(\log n)}+1, ~|B_{i+t}-B_{i}|\leq C'_0a_n\sqrt{n}\right\}$$ and denote $\Delta_n:=\{\max_{s\leq \eta_n}|U(s)-B_s|\leq a_n\sqrt{n}\}.$
Thus conditionally on the event $J_n,$ it is true that
\begin{eqnarray}\label{l01+}
&&\mbfE_{\mu}\Big(\max_{s\leq n}n^{-\frac{1}{2}}|S(s)-(B_s+M(s))|,\Delta_n\cap\Theta_{n}\Big)\no
\leq\mbfE_{\mu}\Big(C_0a_n,\Delta_n\cap\Theta_{n}\Big),\no
\end{eqnarray}where $C_{0}:=1+2C'_{0}$.
Combining with the fact that $f$ is bounded Lipschitz, we have
\begin{eqnarray}\label{l01}
\mbfE_{\mu}\Big(f(S_{(n)})\Big)\no
&\leq&\mbfE_{\mu}\Big(f(B_{(n)}+M_{(n)})+2Ka_n, \Delta_n\cap\Theta_n\Big)+L\mbfP_{\mu}((\Delta_n\cap\Theta_n)^c).\no
\\&\leq&\mbfE_{\mu}\Big(f(B_{(n)}+M_{(n)})\Big)+2Ka_n+L[\mbfP_{\mu}(\Delta^c_n)+\mbfP_{\mu}(\Theta^c_n)]\no
\\&=&\mbfE_{\mu}\Big(f(B+M_{(n)})\Big)+2Ka_n+L[\mbfP_{\mu}(\Delta^c_n)+\mbfP_{\mu}(\Theta^c_n)].
\end{eqnarray}
 For $n$ large enough such that $\frac{na^2_n}{C_1\log n}\geq 1,$ by Theorem \uppercase\expandafter{\romannumeral2} we can see
\begin{eqnarray}
\mbfP_{\mu}\left(\Theta^c_{n}\right)&\leq& \sum_{i=1}^n \mbfP_{\mu}\left(\sup_{|t|\leq (C_1\log n)^{-1}na^2_n+1}|B_{i+t}-B_{i}|> C'_0a_n\sqrt{n}\right)\no
\\&\leq& 2n\mbfP_{\mu}\left(\sup_{t\leq 2(C_1\log n)^{-1}na^2_n}|B_t|> C'_0a_n\sqrt{n}\right)\no
\\&\leq& 2n\cdot D_2e^{-\frac{C'^2_0}{2(C_1\log n)^{-1}D_1}}.\no
\end{eqnarray}
In the rest of the proof, we always let $D_1=3$ (introduced in Theorem \uppercase\expandafter{\romannumeral2}) and write the corresponding $D_2$ as $D$ when we take $D_1=3.$ Hence the above inequality can be rewritten as
\begin{eqnarray}\label{l02+}
\mbfP_{\mu}\left(\Theta^c_{n}\right)\leq 2Dn^{1-\frac{C_1C'^2_0}{6}}.
\end{eqnarray}
Recalling that $C_1C'^2_0>3\alpha$, \eqref{l02+} implies that
$$\varlimsup_{n\ra\infty}\frac{\mbfP_{\mu}\left(\Theta^c_{n}\right)}{n(a_n\sqrt{n}) ^{-\alpha}}=0.$$
Combining with $\eqref{ab}$ and $\eqref{del},$ on event $J_n\cap G_n,$ for $n$ large enough we have
\begin{eqnarray}\label{l02++}\mbfP_{\mu}(\Delta^c_n)+\mbfP_{\mu}(\Theta^c_n)\leq 2C^{(\alpha)}b_n(a_n\sqrt{n})^{-\alpha}.\end{eqnarray}
By \eqref{l01} and \eqref{l02++} we obtain that
\begin{eqnarray}\label{l07}
\mbfE_{\mu}\Big(f(S_{(n)})\Big)\leq \mbfE_{\mu}\Big(f(B+M_{(n)})\Big)+KC_0a_n+2LC^{(\alpha)}b_n(a_n\sqrt{n})^{-\alpha}.\no\end{eqnarray}
Then we show the upper bound in \eqref{lemu} by taking $C_2:=KC_0+2LC^{(\alpha)}$. For the lower bound, we just consider the relationship
\begin{eqnarray}\label{l04}
\mbfE_{\mu}\Big(f(S_{(n)})\Big)
&\geq&\mbfE_{\mu}\Big(f(B_{(n)}+M_{(n)})-2Ka_n, \Delta_n\cap\Theta_n\Big)\no
\\&\geq&\mbfE_{\mu}\Big(f(B_{(n)}+M_{(n)}), \Delta_n\cap\Theta_n\Big)-2Ka_n\mbfP_{\mu}\Big(\Delta_n\cap\Theta_n\Big)\no
\\&\geq&\mbfE_{\mu}\Big(f(B_{(n)}+M_{(n)})\Big)-L[\mbfP_{\mu}(\Delta^c_n)+\mbfP_{\mu}(\Theta^c_n)]-2Ka_n\no
\\&=&\mbfE_{\mu}\Big(f(B+M_{(n)})\Big)-L[\mbfP_{\mu}(\Delta^c_n)+\mbfP_{\mu}(\Theta^c_n)]-2Ka_n.\no
\end{eqnarray}
Combining with \eqref{l02++} we get the lower bound in \eqref{lemu}.\qed

Based on this lemma, we can prove the main theorems.

\noindent{\it Proof of Theorem \ref{t1}.}
~Recall the definition of $\tilde{a}$ in \eqref{a2a}. In this proof we take the $C_1$ in Lemma \ref{L1} satisfy that
\begin{eqnarray}\label{thm00}
C_1\in\left(0, \frac{\tilde{a}^2}{2\sqrt{6(\beta_2-1)}}\right).
\end{eqnarray}
and still choose $C'_0\geq\sqrt{3\alpha/C_1}.$ 

The above lemma means that
\begin{eqnarray}\label{mp1a}F_{n}(x)\leq \mbfP\left(\mbfE_{\mu}\Big(f(B+M_{(n)})\Big)\leq x+C_2y_n\right)+\mbfP((J_nG_n)^c).\end{eqnarray}

Define a random broken line $\hat M(s):=M(\sigma^{-2}s).$ Since $\{M_n\}$ is an i.i.d. random walk, according to Theorem \uppercase\expandafter{\romannumeral1} we can construct a standard Brownian motion $W$ under $\mbfP$ such that
\begin{eqnarray}\label{mp2a}\mbfP\left(\Psi^c_n\right)&\leq& C^{(\beta_1)}
\frac{\sum_{k=1}^n\mbfE(|M_k-M_{k-1}|^{\beta_1})}{(a_n\sqrt{n})^{\alpha}}
=\frac{nC^{(\beta_1)}\mbfE(|M_1|^{\beta_1})}{(a_n\sqrt{n})^{\alpha}},\end{eqnarray}
where $\Psi_n:=\{\sup_{s\leq n\sigma^2}|\hat M(s)-W_s|\leq a_n\sqrt{n}\}.$ Hence on the event $\Psi_n$ (i.e. $\mu\in\Psi_n$), we see
$$\mbfE_{\mu}(f(B+M_{(n)}))-\mbfE_{\mu}(f(B+\hat{W}_{(n)}))\leq Ka_n$$ from the Lipschitz property of $f$, where $\hat{W}_{(n)}(s)=n^{-\frac{1}{2}}W_{\sigma^2s n}.$ Then we have
\begin{eqnarray}\label{mp4a}&&\mbfP\left(\mbfE_{\mu}\Big(f(B+M_{(n)})\Big)\leq x+C_2y_n\right)\no
\\&\leq&\mbfP\left(\mbfE_{\mu}\Big(f(B+\hat{W})\Big)\leq x+C_2y_n+Ka_n\right)\no
\\&=&\mbfP\left(\mbfE_{\mu}\Big(f(B+\sigma W)\Big)\leq x+C_2y_n+Ka_n\right).
\end{eqnarray}
Note that the realization of $W$ totally determines the quantity of $\mbfE_{\mu}\big(f(B+\sigma W)\big)$. That is to say, the impact of $\mu$ on the conditional expectation is only determined by $W.$ Hence in the following we can write this conditional expectation as $\mbfE\big(f(B+\sigma W)|W\big)$.
According to \eqref{mp1a} and \eqref{mp4a}, it is true that
\begin{eqnarray}\label{mp5a}F_{n}(x)\leq \mbfP\left(\mbfE\big(f(B+\sigma W)|W\big)\leq x+Ka_n+C_2y_n\right)+\mbfP(J_n^c)+\mbfP(G_n^c)+\mbfP(\Psi_n^c).~\end{eqnarray}
Now it is turn to estimate $\mbfP(J_n^c)$ and $\mbfP(G_n^c).$ Recall that $\mbfE(\eta_1)=1.$
By Theorem \uppercase\expandafter{\romannumeral1} we can construct a standard Brownian motion $W^{(\eta)}$ such that
$$\mbfP\left(\max_{i\leq n}|(\eta_{i}-i)-W^{(\eta)}_{i}|\geq \frac{na^2_n}{2C_1\log n}\right)\leq C^{(\beta_2)}\mbfE(|\eta_1-1|^{\beta_2})n\left(\frac{2C_1\log n}{na^2_n}\right)^{\beta_2}.$$
Recall that $J_n:=\left\{\forall_{i\leq n}|\eta_{i}-i|\leq \frac{na^2_n}{C_1\log n}\right\}$. Hence by Theorem \uppercase\expandafter{\romannumeral2} we can see
\begin{eqnarray}\label{thm00+}\mbfP(J^c_n)&\leq& \mbfP\left(\max_{i\leq n}|(\eta_{i}-i)-W^{(\eta)}_{i}|\geq \frac{na^2_n}{2C_1\log n}\right)+\mbfP\left(\max_{i\leq n}|W^{(\eta)}_{i}|\geq \frac{na^2_n}{2C_1\log n}\right)\no
\\&\leq&C^{(\beta_2)}\mbfE(|\eta_1-1|^{\beta_2})n\left(\frac{2C_1\log n}{na^2_n}\right)^{\beta_2}+D\exp\left\{-\frac{na^4_n}{12C^2_1(\log n)^2}\right\}.\end{eqnarray}
By the assumption \eqref{ab} and the relationship \eqref{thm00} one can see the second term of the last line in \eqref{thm00+} have a smaller order than the first one,
which means that we can find a constant $C_J$ such that for any $n\in\bfN^+,$ $$\mbfP(J^c_n)\leq C_Jn\left(\frac{\log n}{na^2_n}\right)^{\beta_2}.$$
Note that $\{b_n-n\mbfE(\zeta_1)\}_{n}$ is a centered i.i.d. random walk with finite $\beta_3$-moment. Therefore, with the same method on the estimate of $\mbfP(J^c_n)$ one can see there exists a constant $C_G$ such that
$\mbfP(G^c_n)\leq 4^{-\beta_3}C_Gn\left(\frac{b_n-n\mbfE(\zeta_1)}{2}\right)^{-\beta_3}, \forall n\in\bfN^+.$
 Note that we have set $b_n\geq 2n\mbfE(\zeta_1),$ hence this estimate do not require a relationship like \eqref{thm00} and
 $$\mbfP(G^c_n)\leq C_Gnb^{-\beta_3}_n.$$
Recalling \eqref{mp2a} and \eqref{thm00+} we can find a constant $K_2$ such that
\begin{eqnarray}\label{mp6a}\mbfP(\Psi^c_n)+\mbfP(J^c_n)+\mbfP(G^c_n)\leq K_2r_n, ~~\forall n\in\bfN^+.\end{eqnarray}
Denote $K_1:=C_2+K.$ Combining \eqref{mp5a} with \eqref{mp6a} we see that
\begin{eqnarray}\label{biao1}
F_{n}(x)\leq \mbfP\left(\mbfE\big(f(B+\sigma W)|W\big)\leq x+K_1y_n\right)+K_2r_n,\end{eqnarray}
which completes the upper bound \eqref{thu}.
For the lower bound \eqref{thl}, we should note that
\begin{eqnarray}\label{mp7a}F_{n}(x)
&\geq& \mbfP\left(\mbfE_{\mu}\Big(f(B+M_{(n)})\Big)\leq x-C_2y_n, J_n,G_n\right)\no
\\&\geq& \mbfP\left(\mbfE_{\mu}\Big(f(B+M_{(n)})\Big)\leq x-C_2y_n\right)-\mbfP((J_nG_n)^c).\no\end{eqnarray}
and
\begin{eqnarray}\label{mp7a}&&\mbfP\left(\mbfE_{\mu}\Big(f(B+M_{(n)})\Big)\leq x-C_2y_n, \Psi_n\right)\no
\\&\geq& \mbfP\left(\mbfE\big(f(B+\sigma W)|W\big)\leq x-C_2y_n-Ka_n, \Psi_n\right)\no
\\&\geq& \mbfP\left(\mbfE\big(f(B+\sigma W)|W\big)\leq x-K_1y_n\right)-\mbfP(\Psi_n^c)\no.\end{eqnarray}
Here we omit the details because it is similar to the proof of the upper bound. \qed

 \ack
We would like to thank the referee(s) for the valuable suggestions and useful comments.
This work is supported by the Fundamental Research Funds for the Central Universities (NO.2232021D-30) and the National Natural Science Foundation of China (NO.11971062).

\end{document}